\documentclass[11pt,a4paper]{article}
\usepackage{amsmath}
\usepackage{amsthm}
\usepackage{amssymb}
\usepackage{amsfonts}
\usepackage[colorlinks]{hyperref}

\usepackage{textcomp}
\usepackage{eurosym}
\usepackage[noadjust]{cite}
\usepackage{filecontents}    
\begin{document}
\vspace{3cm}
\begin{center}
\vspace{2cm} {\bf\large On Fixed Points in the Setting of $C^*$-Algebra-Valued Controlled $F_c$-Metric Type Spaces
              }\vspace{1cm}\\
                {\bf G. Kalpana$^{1,*}$ and Z. Sumaiya Tasneem$^1$}\\
                $^1$Department of Mathematics,\\
                SSN College of Engineering, \\
                Kalavakkam, Chennai-603 110, India.\\
               $^*$Correspondence: kalpanag@ssn.edu.in\\
               sumaiyatasneemz@ssn.edu.in\\
\vspace {1cm}
\end {center}

\begin{abstract}
In the present article,\;we first examine the conception of $C^*$-algebra-valued controlled $F_c$-metric type spaces as a generalization of $F$-cone metric spaces over banach algebra. Further, we prove some fixed point theorem with different contractive conditions in the framework of $C^*$-algebra-valued controlled $F_c$-metric type spaces.\;Secondly,\;we furnish an example by means of the acquired result.
 
\vskip1em \noindent \textbf{2010 AMS Classification:} 47H10,\;54H25.

\vskip1em \noindent \textbf{Keywords and Phrases:} $C^\ast$-algebra;\;$C^\ast$-algebra-valued controlled $F_c$-metric type
spaces; Contractive mapping;\;Fixed point theorem.

\vskip1em \noindent \textbf{Article Type:} Research Article. \ \

\end{abstract}
\section{Introduction}
The conception of $b$-metric space was initiated by Bakhtin \cite{Ba89} as a generalization of metric spaces.\;In 1994,\;Matthews \cite{MT94} proposed the concept of partial metric spaces where the self-distance of any point need not be zero.\;Tayyab Kamran et al.\,\cite{KSA17} introduced a new type of metric spaces,\;namely extended b-metric spaces by replacing the constant s by a function $\theta(x,y)$ depending on the parameters of the left-hand side of the triangle inequality.\;Nabil Mlaiki et al.\,\cite{MAAS18} proved banach contraction principle in the setting of controlled metric type spaces which is a generalization of extended $b$-metric space.\;For more engrossing results in extended $b$-metric spaces,\;the reader may refer to \cite{SKP18,SAM18,KSL18,AFK18,PTA19, AFKE18,AKO18}.\;In \cite{MU18},\;Aiman Mukheimer have recently examined the hypothesis of extended partial $S_b$-metric spaces.

On the other hand,\;Fernandez et al. \cite{FMRS17} established the notion of $F$-cone metric space over banach algebra and investigated the existence and uniqueness of the fixed point under the same metric.\;In \cite{ZLH14},\;Ma initiated the concept of $C^*$-algebra-valued metric spaces where the set of real numbers is replaced by the set of all positive elements of a unital $C^*$-algebra.\;For further probes on $C^*$-algebra,\;we refer to \cite{Ba16,BK15,EA18,KK18,KAS19,QLZ16,SD16,KS18,MG90,MZJ15}.

As noted above,\;a vigorous research on fixed point results in $C^*$-algebra-valued metric spaces,\;extended $b$-metric spaces and controlled metric type spaces has been developed in the past few years,\;we focus our study on the concept of $C^*$-algebra-valued controlled $F_c$-metric type spaces in the
present paper and prove fixed point theorem with disparate contractive condition.
\section{Preliminaries}
To start with,\;we recollect some necessary definitions which will be utilized in the main theorem.

Throughout this paper,\;$\mathbb{A}$ denotes an unital $C^*$-algebra.\;Set $\mathbb{A}_h$ = $\{z\in\mathbb{A}:z=z^\ast
\}.$\;We call an element $z \in \mathbb{A}$ a positive element,\;denote it by $\theta_\mathbb{A} \preceq z,$ if $z\in \mathbb{A}_h$ and
$\sigma (z) \subseteq [0,\;\infty)$,\;where $\theta_\mathbb{A}$ is a zero element in $\mathbb{A}$ and $\sigma (z)$ is the spectrum of $z$.\;There is a natural partial ordering on $\mathbb{A}_h$ given by $z \preceq w$ if and only if $\theta_\mathbb{A} \preceq w - z$.\;We denote $\mathbb{A}_{+}$ and $\mathbb{A}^{'}_I$ as $\{z \in \mathbb{A}: \theta_\mathbb{A} \preceq z\}$ and the
set $\{z \in \mathbb A: zw = wz,\;\forall w \in \mathbb{A}\}$ and $|z| = (z^*z)^\frac{1}{2}$ respectively.
\newtheorem{lem}[subsection]{Lemma}
\newtheorem{defn}[subsection]{Definition}
\begin{defn}\cite{FMRS17}
Let $X$ be a nonempty set.\;A function $F: X^3\to A$ is called $F$-cone metric on $X$ if for any $\alpha,\beta,\gamma,\delta \in X$,\;the following conditions hold:\\
$1.$ $\alpha = \beta = \gamma$ if and only if $F(\alpha,\alpha,\alpha) = F(\beta,\beta,\beta) = F(\gamma,\gamma,\gamma) = F(\alpha,\beta,\gamma)$;\\
$2.$ $\theta \preceq F(\alpha,\alpha,\alpha) \preceq F(\alpha,\alpha,\beta) \preceq F(\alpha,\beta,\gamma)$ for all $\alpha,\beta,\gamma \in X$ with $\alpha\neq\beta\neq \gamma$;\\
$3.$ $F(\alpha,\beta,\gamma) \preceq s[F(\alpha,\alpha,\delta)+ F(\beta,\beta,\delta) + F(\gamma,\gamma,\delta)] - F(\delta,\delta,\delta)$.\\
Then the pair $(X,F)$ is called an $F$-cone metric space over Banach Algebra $A$.\;The number $s\geq 1$ is called the coefficient of $(X,F)$.
\end{defn}
\begin{defn}\cite{KK18}
Let $X$ be a nonempty set and $A \in \mathbb{A}^{'}$ such
that $A \succeq I_\mathbb{A}$.\;Suppose the mapping $S_b : X \times X
\times X \to \mathbb{A}$ satisfies:\\
$1.$ $\theta_\mathbb{A} \preceq S_b(\alpha,\beta,\gamma)$ for all $\alpha,\beta,\gamma \in X$ with $\alpha \neq \beta \neq \gamma \neq \alpha;$\\
$2.$ $S_b(\alpha,\beta,\gamma) = \theta_\mathbb{A}$ if and only if $ \alpha = \beta = \gamma;$\\
$3.$ $S_b(\alpha,\beta,\gamma) \preceq A [S_b(\alpha,\alpha,\delta)+ S_b(\beta,\beta,\delta)+ S_b(\gamma,\gamma,\delta)]$ for all $\alpha,\beta,\gamma,\delta \in X.$
\end{defn}
Then $S_b$ is said to be $C^*$-algebra-valued $S_b$-metric on $X$
and $(X, \mathbb{A},S_b)$ is said to be a $C^\ast$-algebra-valued
$S_b$-metric space.

\begin{defn}\cite{MAAS18}
Given a non-empty set $X$ and $\delta:X\times X\to [1,\infty)$.\;A function $d:X\times X\to[0,\infty)$ is called a controlled metric type if:\\[3pt]
$1.$ $d(\alpha,\beta) = 0$ if and only if $\alpha = \beta$;\\[3pt]
$2.$ $d(\alpha,\beta) = d(\beta,\alpha)$;\\[3pt]
$3.$ $d(\alpha,\beta) \leq \delta(\alpha,\gamma)d(\alpha,\gamma) + \delta(\gamma,\beta)d(\gamma,\beta)$,\;for all $\alpha,\beta,\gamma\in X$.\\[3pt]
The pair $(X,d)$ is called a controlled metric type space.
\end{defn}

\section{Main Results}
In this main segment,\;as a generalization of $F$-cone metric space over Banach algebra,\;we introduce the notion of $C^*$-algebra valued controlled $F_c$-metric type spaces and furnish an example of the underlying spaces.\\
Hereinafter,\;$\mathbb{A}_I'$ will denote the set $\{z \in \mathbb A: zw = wz,\;\forall w \in \mathbb{A}\;\text{and}\; z\succeq I_\mathbb{A}\}$ respectively.
\begin{defn}
Let $X$ be a nonempty set and $C:X \times X \times X \to \mathbb{A}_I'$.\;Suppose the mapping $F_c:X \times X
\times X \to \mathbb{A}$ satisfies:\\
$1.$ $\varpi = \bar{\nu} = \bar{\varsigma}$ if and only if $F_c(\varpi,\varpi,\varpi) = F_c(\bar{\nu},\bar{\nu},\bar{\nu}) = F_c(\bar{\varsigma},\bar{\varsigma},\bar{\varsigma}) = F_c(\varpi,\bar{\nu},\bar{\varsigma})$;\\
$2.$ $\theta_\mathbb{A} \preceq F_c(\varpi,\varpi,\varpi) \preceq F_c(\varpi,\varpi,\bar{\nu})\preceq F_c(\varpi,\bar{\nu},\bar{\varsigma})$;\\
$3.$ $F_c(\varpi,\bar{\nu},\bar{\varsigma}) \preceq C(\varpi,\varpi,\bar{\alpha})F_c(\varpi,\varpi,\bar{\alpha})+ C(\bar{\nu},\bar{\nu},\bar{\alpha})F_c(\bar{\nu},\bar{\nu},\bar{\alpha})+ C(\bar{\varsigma},\bar{\varsigma},\bar{\alpha})F_c(\bar{\varsigma},\bar{\varsigma},\bar{\alpha}) - F_c(\bar{\alpha},\bar{\alpha},\bar{\alpha})$ for all $\varpi,\bar{\nu},\bar{\varsigma},\bar{\alpha} \in X.$
\end{defn}
Then $F_c$ is called a $C^*$-algebra-valued controlled $F_c$- metric type on $X$ and $(X, \mathbb{A},F_c)$ is a $C^\ast$-algebra-valued controlled $F_c$-metric type spaces.

\newtheorem{rem}[subsection]{Remark}
\begin{rem} 
If $C(\varpi,\varpi,\bar{\alpha}) = C(\bar{\nu},\bar{\nu},\bar{\alpha}) = C(\bar{\varsigma},\bar{\varsigma},\bar{\alpha}) = C(\varpi,\bar{\nu},\bar{\varsigma})$ for all $\varpi,\bar{\nu},\bar{\varsigma},\bar{\alpha} \in X$,\;then we get
$$ F_c(\varpi,\bar{\nu},\bar{\varsigma}) \preceq C(\varpi,\bar{\nu},\bar{\varsigma})[F_c(\varpi,\varpi,\bar{\alpha})+ F_c(\bar{\nu},\bar{\nu},\bar{\alpha})+ F_c(\bar{\varsigma},\bar{\varsigma},\bar{\alpha})] - F_c(\bar{\alpha},\bar{\alpha},\bar{\alpha}).$$
In this case,\;$F_c$ is called a $C^*$-algebra-valued extended $F_c$-metric on $X$ and ($X,\mathbb{A},c$) is called a $C^*$-algebra-valued extended $F_c$-metric space.
\end{rem}
\begin{rem}
In a $C^*$-algebra-valued controlled $F_c$-metric type space $(X,\mathbb{A},F_c)$,\;if $\varpi,\bar{\nu},\bar{\varsigma}\in X$ and $F_c(\varpi,\bar{\nu},\bar{\varsigma}) = \theta$,\;then $\varpi =\bar{\nu}=\bar{\varsigma}$,\;but the converse need not be true.
\end{rem}
\begin{defn}
A $C^*$-algebra-valued controlled $F_c$-metric type space $(X,\mathbb{A},F_c)$ is said to be symmetric if it satisfies, $$F_c(\varpi,\varpi,\bar{\nu}) = F_c(\bar{\nu},\bar{\nu},\varpi),\;\text{for all}\;\varpi,\bar{\nu}\in X.$$
\end{defn}
\newtheorem{ex}[subsection]{Example}
\begin{ex}
Let $X=\{0,1,2,\ldots\}$ and $\mathbb{A} = \mathbb{R}^2$.\;If $\alpha,\beta\in \mathbb{A}$ with $\varpi=(\varpi_1,\varpi_2),\,\bar{\nu}=(\bar{\nu}_1,\bar{\nu}_2)$,\;then the addition,\;multipilcation and scalar multipilcation can be defined as follows: 
$$ \varpi + \bar{\nu} = (\varpi_1+\bar{\nu}_1,\;\varpi_2+\bar{\nu}_2),\;k\varpi = (k\varpi_1,k\varpi_2),\;\varpi\bar{\nu} = (\varpi_1\bar{\nu}_1,\varpi_2\bar{\nu}_2) $$
Now define the metric $F_c:X \times X \times X \to \mathbb{A}$ and the control function $C:X \times X \times X \to \mathbb{A}'_I$ as:
$$F_c(\varpi,\bar{\nu},\bar{\varsigma}) = \Bigl(\frac{1}{2}(|\varpi+\bar{\varsigma}|^2 + |\bar{\nu}+\bar{\varsigma}|^2),\;\frac{1}{2}(|\varpi+\bar{\varsigma}|^2 + |\bar{\nu}+\bar{\varsigma}|^2)\Bigr)$$
and
$$C(\varpi,\bar{\nu},\bar{\varsigma}) = \Bigl(|\varpi+\bar{\nu}-\bar{\varsigma}+1|,|\varpi+\bar{\nu}-\bar{\varsigma}+1|\Bigr).$$
It is easy to verify that $F_c$ is a $C^*$-algebra-valued controlled $F_c$-metric type space.\;Indeed for $\varpi =1,\;\bar{\nu} = 2,\;\bar{\varsigma} =3 \;\text{and}\;\bar{\alpha} = 0$,\;we have
\begin{equation*}
\begin{aligned}
F_c(1,2,3) = (20.5,\;20.5) &\succ (1,\;1) [(1,\;1) + (4,\;4) + (9,\;9)]  - (0,0) = (14,\;14)\\[2pt]
&= C(1,2,3)[F_c(1,1,0) + F_c(2,2,0) + F_c(3,3,0)] - F_c(0,0,0).
\end{aligned}
\end{equation*}
Hence $F_c$ is not a $C^*$-algebra-valued extended $F_c$-metric  space.
\end{ex}
\begin{defn}
A sequence $\{\varpi_n\}$ in a $C^*$-algebra-valued controlled $F_c$-metric type space is said to be:\\
(i) convergent sequence $\Longleftrightarrow \exists\; \varpi \in X$ such that $F_c(\varpi_n,\varpi_n,\varpi)\to \theta_\mathbb{A}$ as $n \to \infty$ and we denote it by $\lim\limits_{n \to \infty} \varpi_n = \varpi$;\\
(ii) Cauchy sequence $\Longleftrightarrow F_c(\varpi_n,\varpi_n,\varpi_m)\to \theta_\mathbb{A}$ as $n,m \to \infty$.
\end{defn}
\begin{defn}
A $C^*$-algebra-valued controlled $F_c$-metric type space $(X,\mathbb{A},F_c)$ is said to be complete if every Cauchy sequence is convergent in $X$ with respect to $\mathbb{A}$.
\end{defn}
\newtheorem{thm}[subsection]{Theorem}
\begin{thm}\label{thm1}
Let $(X,\mathbb{A},F_c)$ be a complete symmetric $C^*$-algebra-valued controlled $F_c$-metric type space and suppose $T:X\to X$ is a mapping satisfying the following condition:
\begin{equation}\label{eq1.}
F_c(T\varpi,T\varpi,T\bar{\nu}) \preceq P^*\,F_c(\varpi,\varpi,\bar{\nu})\,P + Q^*\,F_c(\varpi,\varpi,T\varpi)\,Q + R^*\,F_c(\bar{\nu},\bar{\nu},T\bar{\nu})\,R, \;\forall \varpi,\,\bar{\nu},\in X,
\end{equation}
where $P,\,Q,\,R\in \mathbb{A}$ with $\|P\|,\;\|Q\|,\;\|R\|\geq0$ satisfying $\|P\|^2 + \|Q\|^2 + \|R\|^2 < 1$ and for $\varpi_0\in X$,\;choose $\varpi_n=T^n \varpi_0$ assume that
\begin{equation}\label{eq2.}
\sup\limits_{m\geq1}\lim\limits_{i\to\infty} \|C(\varpi_{i+1},\varpi_{i+1},\varpi_{i+2}) C(\varpi_{i+1},\varpi_{i+1},\varpi_m)\| < \frac{1 - \|R\|^2}{\|P\|^2 + \|Q\|^2}.
\end{equation}
In addition,\;for each $\varpi\in X$,\;suppose that 
\begin{equation}\label{eq3}
\lim\limits_{n\to\infty} \|C(\varpi,\varpi,\varpi_n)\| \;\text{and}\;\lim\limits_{n\to\infty} \|C(\varpi_n,\varpi_n,\varpi)\|
\end{equation}
exist and are finite.\;Then $T$ has a unique fixed point in $X$.
\begin{proof}
Let $\varpi_0 \in X$ be arbitrary and define the iterative sequence $\{\varpi_n\}$ by:
\begin{equation}\label{eq4}
\varpi_{n+1} = T\varpi_n = \ldots = T^{n+1}\varpi_0,\;n= 1,2,\ldots. 
\end{equation}
If follows from (\ref{eq1.}) and (\ref{eq4}) that
\begin{equation}
\begin{aligned}
F_c(\varpi_{n},\varpi_{n},\varpi_{n+1}) &= F_c(T\varpi_{n-1},T\varpi_{n-1},T\varpi_{n})\\[2pt]
&\preceq P^* F_c(\varpi_{n-1},\varpi_{n-1},\varpi_{n})\,P + Q^*\,F_c(\varpi_{n-1},\varpi_{n-1},T\varpi_{n-1})\,Q +\\[2pt] &\quad\, R^*\,F_c(\varpi_{n},\varpi_{n},T{\varpi_n})\,R\\[2pt]
\Longleftrightarrow \|F_c(\varpi_{n},\varpi_{n},\varpi_{n+1})\| &\leq \|P^* F_c(\varpi_{n-1},\varpi_{n-1},\varpi_{n})\,P + Q^*\,F_c(\varpi_{n-1},\varpi_{n-1},T\varpi_{n-1})\,Q +\\[2pt] &\quad\, R^*\,F_c(\varpi_{n},\varpi_{n},T{\varpi_n})\,R\|\\[2pt]
&\leq \|P^* F_c(\varpi_{n-1},\varpi_{n-1},\varpi_{n})\,P\| + \|Q^*\,F_c(\varpi_{n-1},\varpi_{n-1},T\varpi_{n-1})\,Q\| +\\[2pt] &\quad\, \|R^*\,F_c(\varpi_{n},\varpi_{n},T{\varpi_n})\,R\|\\[2pt]
&= (\|P\|^2 + \|Q\|^2) \|F_c(\varpi_{n-1},\varpi_{n-1},\varpi_{n})\| + \|R\|^2 \|F_c(\varpi_{n},\varpi_{n},\varpi_{n+1})\|\\[8pt]
\therefore \|F_c(\varpi_{n},\varpi_{n},\varpi_{n+1})\| &\leq \frac{\|P\|^2 + \|Q\|^2}{1-\|R\|^2} \|F_c(\varpi_{n-1},\varpi_{n-1},\varpi_{n})\|
\end{aligned}
\end{equation}
Accordingly we get
\begin{equation}
\begin{aligned}
\|F_c(\varpi_{n},\varpi_{n},\varpi_{n+1})\| &\leq \|S\|^2 \|F_c(\varpi_{n-1},\varpi_{n-1},\varpi_{n})\|\\[5pt]
&= \|S^*S\|\|F_c(\varpi_{n-1},\varpi_{n-1},\varpi_{n})\|\\[5pt]
&\leq \|S^*\|\|F_c(\varpi_{n-1},\varpi_{n-1},\varpi_{n})\|\|S\|\\[5pt]
\Longleftrightarrow F_c(\varpi_{n},\varpi_{n},\varpi_{n+1}) &\preceq S^* F_c(\varpi_{n-1},\varpi_{n-1},\varpi_{n}) S, 
\end{aligned}
\end{equation}
where $\|S\|^2 = \frac{\|P\|^2 + \|Q\|^2}{1-\|R\|^2} < 1$.\;Recursively,\;we find that
\begin{equation}
F_c(\varpi_{n},\varpi_{n},\varpi_{n+1}) \preceq (S^*)^n F_c(\varpi_{n-1},\varpi_{n-1},\varpi_{n}) S^n
\end{equation}
For any $n\geq 1$ and $q\geq 1$,\;we have
\begin{equation*}
\begin{aligned}
F_c(\varpi_{n},\varpi_{n},\varpi_{n+q}) &\preceq C(\varpi_{n},\varpi_{n},\varpi_{n+1}) F_c(\varpi_{n},\varpi_{n},\varpi_{n+1}) + 
C(\varpi_{n},\varpi_{n},\varpi_{n+1}) F_c(\varpi_{n},\varpi_{n},\varpi_{n+1}) + \\[3pt]
&\quad\, C(\varpi_{n+q},\varpi_{n+q},\varpi_{n+1}) F_c(\varpi_{n+q},\varpi_{n+q},\varpi_{n+1}) - F_c(\varpi_{n+1},\varpi_{n+1},\varpi_{n+1})\\[3pt]
&\preceq 2C(\varpi_{n},\varpi_{n},\varpi_{n+1}) F_c(\varpi_{n},\varpi_{n},\varpi_{n+1}) + C(\varpi_{n+q},\varpi_{n+q},\varpi_{n+1}) \\[3pt]
&\quad\, F_c(\varpi_{n+1},\varpi_{n+1},\varpi_{n+q})\\[3pt]
&\preceq 2C(\varpi_{n},\varpi_{n},\varpi_{n+1}) F_c(\varpi_{n},\varpi_{n},\varpi_{n+1}) + C(\varpi_{n+q},\varpi_{n+q},\varpi_{n+1}) \\[3pt]
&\quad\, \bigl[2\,C(\varpi_{n+1},\varpi_{n+1},\varpi_{n+2}) F_c(\varpi_{n+1},\varpi_{n+1},\varpi_{n+2}) + C(\varpi_{n+q},\varpi_{n+q},\varpi_{n+2}) \\[3pt]
&\quad\, F_c(\varpi_{n+2},\varpi_{n+2},\varpi_{n+q})\bigr] - F_c(\varpi_{n+2},\varpi_{n+2},\varpi_{n+2})\\[3pt]
&\quad\, \vdots\\[3pt]
&= 2C(\varpi_{n},\varpi_{n},\varpi_{n+1}) F_c(\varpi_{n},\varpi_{n},\varpi_{n+1}) + \\[3pt]
&\quad\, 2\sum_{i=n+1}^{n+q-2} C(\varpi_{i},\varpi_{i},\varpi_{i+1}) F_c(\varpi_{i},\varpi_{i},\varpi_{i+1})\prod_{j=n+1}^i C(\varpi_{n+q},\varpi_{n+q},\varpi_{j}) + \\[3pt]
&\quad\, \prod_{i=n+1}^{n+q-1} C(\varpi_{n+q},\varpi_{n+q},\varpi_{i}) F_c(\varpi_{n+q-1},\varpi_{n+q-1},\varpi_{n+q}) \\[3pt]
&\preceq 2C(\varpi_{n},\varpi_{n},\varpi_{n+1}) F_c(\varpi_{n},\varpi_{n},\varpi_{n+1}) + \\[3pt]
&\quad\, 2\sum_{i=n+1}^{n+q-1} C(\varpi_{i},\varpi_{i},\varpi_{i+1}) F_c(\varpi_{i},\varpi_{i},\varpi_{i+1})\prod_{j=n+1}^i C(\varpi_{n+q},\varpi_{n+q},\varpi_{j})\\[3pt]
&\preceq 2C(\varpi_{n},\varpi_{n},\varpi_{n+1}) (S^*)^n S_0 S^n + \\[3pt]
&\quad\, 2\sum_{i=n+1}^{n+q-1} C(\varpi_{i},\varpi_{i},\varpi_{i+1}) (S^*)^i S_0 S^i \prod_{j=1}^i C(\varpi_{n+q},\varpi_{n+q},\varpi_{j})\\[3pt]
\end{aligned}
\end{equation*}
\begin{equation*}
\begin{aligned}
&= 2\Bigl(S_0^\frac{1}{2} C(\varpi_{n},\varpi_{n},\varpi_{n+1}) ^\frac{1}{2} S^n\Bigr)^* (S_0^\frac{1}{2} C(\varpi_{n},\varpi_{n},\varpi_{n+1}) ^\frac{1}{2} S^n\Bigr) + \\[2pt]
&\quad\, 2\sum_{i=n+1}^{n+q-1} \Bigl(S_0^\frac{1}{2} \bigl[C(\varpi_{i},\varpi_{i},\varpi_{i+1})\prod_{j=1}^i C(\varpi_{n+q},\varpi_{n+q},\varpi_{j})\bigr]^\frac{1}{2} S^i\Bigr)^* \\[2pt]
&\quad\qquad\quad\; \Bigl(S_0^\frac{1}{2} \bigl[C(\varpi_{i},\varpi_{i},\varpi_{i+1})\prod_{j=1}^i C(\varpi_{n+q},\varpi_{n+q},\varpi_{j})\bigr]^\frac{1}{2} S^i\Bigr)\\[2pt]
&= 2\,|S_0^\frac{1}{2} C(\varpi_{n},\varpi_{n},\varpi_{n+1}) ^\frac{1}{2} S^n|^2 \,+ \\[2pt] 
&\quad\,\,\,2\sum_{i=n+1}^{n+q-1} \Bigl|S_0^\frac{1}{2} \bigl[C(\varpi_{i},\varpi_{i},\varpi_{i+1})\prod_{j=1}^i C(\varpi_{n+q},\varpi_{n+q},\varpi_{j})\bigr]^\frac{1}{2} S^i\Bigr|^2\\[2pt]
&\preceq 2 \|S_0\| \Bigl[\|C(\varpi_{n},\varpi_{n},\varpi_{n+1})\|\, \|S\|^{2n} I_\mathbb{A} + \\[2pt]
&\quad\;\; \|C(\varpi_{i},\varpi_{i},\varpi_{i+1})\prod_{j=1}^i C(\varpi_{n+q},\varpi_{n+q},\varpi_{j})\| \|S\|^{2i} I_\mathbb{A} \Bigr]
\end{aligned}
\end{equation*}
where $I_\mathbb{A}$ is the unit element in $\mathbb{A}$ and $c(\varpi_1,\varpi_1,\varpi_0) = S_0$ for some $S_0 \in \mathbb{A}$. Let $Y_m = \sum_{i=1}^m \|S\|^{2i}\,\|C(\varpi_{i},\varpi_{i},\varpi_{i+1})\prod\limits_{j=1}^i C(\varpi_{n+q},\varpi_{n+q},\varpi_{j})\|$.\;Consequently the above inequality implies, 
\begin{equation}
F_c(\varpi_{n},\varpi_{n},\varpi_{n+q}) \leq 2\|S_0\| \Bigl[\|C(\varpi_{n},\varpi_{n},\varpi_{n+1})\|\, \|S\|^{2n} + (Y_{n+q-1} - Y_n)\Bigr] I_\mathbb{A}
\end{equation}
The ratio test jointly with (\ref{eq2.}) implies that the limit of the sequence $\{Y_n\}$ exists and so $\{Y_n\}$ is Cauchy.\;Letting $n\to\infty$ in the inequality above,\;we get
\begin{equation}
\lim\limits_{n\to\infty} F_c(\varpi_{n},\varpi_{n},\varpi_{n+q}) = \theta_\mathbb{A}.
\end{equation}
Wherefore the sequence $\{\varpi_n\}$ is Cauchy with respect to $\mathbb{A}$.\;Since $(X,\mathbb{A},F_c)$ is a complete $C^*$-algebra-valued controlled $F_c$-metric type space,\;there exists a point $\varpi\in X$ such that 
\begin{equation}\label{eq11}
\lim\limits_{n\to\infty} F_c(\varpi_{n},\varpi_{n},\varpi) = \theta_\mathbb{A}.
\end{equation}
Consider,
\begin{equation*}
\begin{aligned}
F_c(\varpi,\varpi,\varpi_{n+1}) &\preceq 2C(\varpi,\varpi,\varpi_{n})F_c(\varpi,\varpi,\varpi_{n}) + C(\varpi_{n+1},\varpi_{n+1},\varpi_{n})F_c(\varpi_n,\varpi_n,\varpi_{n+1})\\[2pt]
&\quad\, - F_c(\varpi_n,\varpi_n,\varpi_{n}) \\[2pt]
\Longleftrightarrow
\|F_c(\varpi,\varpi,\varpi_{n+1})\| &\leq 2\|C(\varpi,\varpi,\varpi_{n})\| \|F_c(\varpi,\varpi,\varpi_{n})\| + \|C(\varpi_{n+1},\varpi_{n+1},\varpi_{n})\| \\[2pt]
&\quad\,\;\|F_c(\varpi_n,\varpi_n,\varpi_{n+1})\|
\end{aligned}
\end{equation*}
It yields from (\ref{eq3.}) and (\ref{eq11}) that
\begin{equation}\label{eq12}
\lim\limits_{n\to\infty} \|F_c(\varpi,\varpi,\varpi_{n+1})\| = 0.
\end{equation}
Hence
\begin{equation*}
\begin{aligned}
\|F_c(\varpi,\varpi,T\varpi)\| &\leq 2\|C(\varpi,\varpi,\varpi_{n+1})\| \|F_c(\varpi,\varpi,\varpi_{n+1})\| + \|C(T\varpi,T\varpi,\varpi_{n+1})\| \\[2pt]
&\quad\,\;\|F_c(\varpi_{n+1},\varpi_{n+1},T\varpi)\|\\[2pt]
&= 2\|C(\varpi,\varpi,\varpi_{n+1})\| \|F_c(\varpi,\varpi,\varpi_{n+1})\| + \|C(T\varpi,T\varpi,\varpi_{n+1})\| \\[2pt]  
&\quad\,\; \|F_c(T^{n+1}\varpi,T^{n+1}\varpi,T\varpi)\|
\end{aligned}
\end{equation*}
Regarding (\ref{eq12}),\;we get $\|F_c(\varpi,\varpi,\varpi_{n+1})\| \to 0$ as $n\to\infty$.\;Since $T^n\to x$ and from continuity of $T$,\;we acquire $T^{n+1}\to Tx$ i.e.,\;$\|F_c(T^{n+1}\varpi,T^{n+1}\varpi,T\varpi)\| \to 0,$\;as $n\to\infty$.\;Thus
$$ \lim\limits_{n\to\infty} \|F_c(\varpi,\varpi,T\varpi)\| = 0$$
$$\Longleftrightarrow \lim\limits_{n\to\infty} F_c(\varpi,\varpi,T\varpi) = \theta_\mathbb{A}.$$ 
Hence $T\varpi = \varpi$ i.e., $\varpi$ is a fixed point of $T$.\;Now to prove uniqueness,\;let $\bar{\nu}\neq \varpi$ be another fixed point of $T$.\;Taking the expression (\ref{eq1.}) into account,\;we have
\begin{equation*}
\begin{aligned}
F_c(\varpi,\varpi,\bar{\nu}) &= F_c(T\varpi,T\varpi,T\bar{\nu}) \\[2pt]
&\preceq P^*\,F_c(\varpi,\varpi,\bar{\nu})\,P + Q^*\,F_c(\varpi,\varpi,T\varpi)\,Q + R^*\,F_c(\bar{\nu},\bar{\nu},T\bar{\nu})\,R\\[2pt]
&= P^*\,F_c(\varpi,\varpi,\bar{\nu})\,P + Q^*\,F_c(\varpi,\varpi,\varpi)\,Q + R^*\,F_c(\bar{\nu},\bar{\nu},\bar{\nu})\,R\\[2pt]
&\preceq P^*\,F_c(\varpi,\varpi,\bar{\nu})\,P + Q^*\,F_c(\varpi,\varpi,\bar{\nu})\,Q + R^*\,F_c(\bar{\nu},\bar{\nu},\varpi)\,R \\[2pt]
\|F_c(\varpi,\varpi,\bar{\nu})\| &\leq (\|P\|^2 + \|Q\|^2)F_c(\varpi,\varpi,\bar{\nu}) + \|R\|^2 \|F_c(\bar{\nu},\bar{\nu},\varpi)\|\\[2pt]
\|F_c(\varpi,\varpi,\bar{\nu})\| &\leq \frac{\|R\|^2}{(1 - \|P\|^2 - \|Q\|^2)} \|F_c(\bar{\nu},\bar{\nu},\varpi)\|\\[2pt]
&< \|F_c(\bar{\nu},\bar{\nu},\varpi)\| = \|F_c(\varpi,\varpi,\bar{\nu})\|
\end{aligned}
\end{equation*}
which is a contradiction.\;Hence the fixed point is unique.
\end{proof}
\end{thm}
In Theorem (\ref{thm1}),\;if we take $Q = R = \theta$,\;then the above theorem reduces to a Banach contraction principle,\;which can be stated as follows:
\newtheorem{cor}[subsection]{Corollary}
\begin{cor}\label{cor}
Let $(X,\mathbb{A},F_c)$ be a complete $C^*$-algebra-valued controlled $F_c$-metric type space and suppose $T:X\to X$ is a mapping satisfying the following condition:
\begin{equation}\label{eq1}
F_c(T\varpi,T\varpi,T\bar{\nu}) \preceq P^*\,F_c(\varpi,\varpi,\bar{\nu})\,P,\;\;\forall \varpi,\;\bar{\nu},\in X,
\end{equation}
where $P\in \mathbb{A}$ with $0\leq\|P\|<1$  and for $\varpi_0\in X$,\;choose $\varpi_n=T^n \varpi_0$ assume that
\begin{equation}\label{eq2}
\sup\limits_{m\geq1}\lim\limits_{i\to\infty} \|C(\varpi_{i+1},\varpi_{i+1},\varpi_{i+2}) C(\varpi_{i+1},\varpi_{i+1},\varpi_m)\| < \frac{1}{\|P\|^2}.
\end{equation}
In addition,\;for each $\varpi\in X$,\;suppose that 
\begin{equation}\label{eq3.}
\lim\limits_{n\to\infty} \|C(\varpi,\varpi,\varpi_n)\| \;\text{and}\;\lim\limits_{n\to\infty} \|C(\varpi_n,\varpi_n,\varpi)\|
\end{equation}
exist and are finite.\;Then $T$ has a unique fixed point in $X$.
\end{cor}
\begin{ex}
Let $X = [0,4]$ and $\mathbb{A} = M_2(\mathbb{R})$ be the set of all $2\times2$ matrices under usual addition,\;multiplication and scalar multiplication.\;Define $F_c: X\times X\times X \to \mathbb{A}$ as follows:
\begin{equation*}
F_c(\varpi,\bar{\nu},\bar{\varsigma}) = 
\begin{pmatrix}
\text{max}\,\{\varpi,\bar{\varsigma}\} + \text{max}\,\{\bar{\nu},\bar{\varsigma}\} & 0 \\
0 & \text{max}\,\{\varpi,\bar{\varsigma}\} + \text{max}\,\{\bar{\nu},\bar{\varsigma}\}
\end{pmatrix}
\end{equation*}
Hence $(X,\mathbb{A},F_c)$ is a $C^*$-algebra-valued controlled $F_c$-metric type space with $C(\varpi,\bar{\nu},\bar{\varsigma}) = 2 + \text{max}\,\{\varpi,\bar{\nu},\bar{\varsigma}\}$. Now for any $A\in \mathbb{A}$,\;we define its norm as $\|A\| = \max\limits_{1\leq i\leq4}\{|a_i|\}$.\;Let $T: X\to X$ be defined as $T\varpi = \frac{\varpi}{8}$.\;Then 
\begin{equation*}
\begin{aligned}
F_c(T\varpi,T\varpi,T\bar{\nu}) &= F_c(\frac{\varpi}{8},\frac{\varpi}{8},\frac{\bar{\nu}}{8}) \\[2pt]
&= \begin{pmatrix}2\text{max}\,\{\frac{\varpi}{8},\frac{\bar{\nu}}{8}\} & 0 \\
0 & 2\text{max}\,\{\frac{\varpi}{8},\frac{\bar{\nu}}{8}\} \end{pmatrix}\\[2pt]
&= P^* F_c(\varpi,\varpi,\bar{\nu}) P
\end{aligned}
\end{equation*}
where $P = \begin{pmatrix} \frac{1}{2\sqrt{2}} & 0 \\ 0 & \frac{1}{2\sqrt{2}} \end{pmatrix}$ with $\|P\| = \frac{1}{2\sqrt{2}} < 1.$\;Now consider
\begin{equation*}
\begin{aligned}
C(\varpi_{i+1},\varpi_{i+1},\varpi_{i+2}) &= C(T^{i+1}\varpi,T^{i+1}\varpi,T^{i+2}\varpi)\\[2pt]
&= C(\frac{\varpi}{8^{i+1}},\frac{\varpi}{8^{i+1}},\frac{\varpi}{8^{i+2}}) \\[2pt]
&= \begin{pmatrix} 2 + \text{max}\,\{\frac{\varpi}{8^{i+1}},\frac{\varpi}{8^{i+1}},\frac{\varpi}{8^{i+2}}\} & 0 \\ 0 & 2 + \text{max}\,\{\frac{\varpi}{8^{i+1}},\frac{\varpi}{8^{i+1}},\frac{\varpi}{8^{i+2}}\}  \end{pmatrix}
\end{aligned}
\end{equation*}
Similarly,
$$ C(\varpi_{i+1},\varpi_{i+1},\varpi_{m}) = \begin{pmatrix} 2 + \text{max}\,\{\frac{\varpi}{8^{i+1}},\frac{\varpi}{8^{i+1}},\frac{\varpi}{8^{m}}\} & 0 \\ 0 & 2 + \text{max}\,\{\frac{\varpi}{8^{i+1}},\frac{\varpi}{8^{i+1}},\frac{\varpi}{8^{m}}\}  \end{pmatrix}
$$
Thus
\begin{equation*}
\begin{aligned}
\lim\limits_{i\to\infty} & \|C(\varpi_{i+1},\varpi_{i+1},\varpi_{i+2})\;C(\varpi_{i+1},\varpi_{i+1},\varpi_{m})\| 
 \\[2pt]
&= \lim\limits_{i\to\infty}\Bigl\|\begin{pmatrix} (2 + \frac{\varpi}{8^{i+1}})(2 + \text{max}\,(\frac{\varpi}{8^{i+1}},\frac{\varpi}{8^m}) & 0 \\ 0 & (2 + \frac{\varpi}{8^{i+1}})(2 + max(\frac{\varpi}{8^{i+1}},\frac{\varpi}{8^m}) \end{pmatrix}\Bigr\| \\[2pt]
&= \lim\limits_{i\to\infty} (2 + \frac{\varpi}{8^{i+1}})(2 + \text{max}\,(\frac{\varpi}{8^{i+1}},\frac{\varpi}{8^m}) = 4 + \frac{2\varpi}{8^m}
\end{aligned}
\end{equation*}
and $$\sup_{m\geq 1} \lim\limits_{i\to\infty}\|C(\varpi_{i+1},v_{i+1},\varpi_{i+2})\;C(\varpi_{i+1},\varpi_{i+1},\varpi_{m})\| = 4 + \frac{2\varpi}{8} < 8 = \frac{1}{\|P\|^2}.$$ 
Thus $T$ satisfies all the conditions of Theorem (\ref{thm1}),\;hence it has a unique fixed point which is $\varpi = 0.$
\end{ex}

\section{Conclusion}
In this manuscript,\;we have analyzed the structure of $C^*$-algebra-valued controlled $F_c$-metric type spaces and acquired some fixed point theorem under different contractive conditions of the underlying spaces.\;Further,\;an example is conferred to show the effectiveness of the established result.



\begin{thebibliography}{label}
\bibitem{FMRS17}{\rm Fernandez, J., Malviya, N., Radenovi\'{c}, S., Saxena, K.: $F$-cone metric spaces over banach algebra: Fixed Point Theory Appl. 2017, 7 (2017).}
\bibitem{MU18}{\rm Mukheimer, A.:\;Extended partial $S_b$-metric spaces,\;Axioms 2018,\;7,\;87.}
\bibitem{KSL18}{\rm Karapinar, E., Sumati Kumari, P., Lateef, D.:\;A New Approach to the Solution of the Fredholm Integral Equation via a Fixed Point on Extended $b$-Metric Spaces. Symmetry 2018, 10, 512.}
\bibitem{PTA19}{\rm Panda, S.K, Tassaddiq, A., Agarwal, R.P.:\;A New Approach to the Solution of the Non-Linear Integral Equations via Various $\mathbb{F}_{B_e}$-Contractions. Symmetry 2019, 11, 206.}
\bibitem{AFKE18}{\rm Alqahtani, B., Fulga, A., Karapinar, E.:\;Common fixed point results on extended $b$-metric space: J. Inequal. Appl. 2018, 2018, 158.}
\bibitem{AKO18}{\rm Alqahtani, B., Karapinar, E., Ozturk, A.:\;On $(\alpha,\psi)-K$-contractions in the extended $b$-metric space: Filomat 2018, 32, 15.}
\bibitem{AFK18}{\rm Alqahtani, B., Fulga, A., Karapinar, E.:\;Non-Unique Fixed Point Results in Extended $b$-Metric Space. Mathematics 2018, 6, 68.}
\bibitem{SKP18}{\rm Samreen, M., Kamran, T., Postolache, M.:\;Extended $b$-metric space, extended $b$-comparison function and nonlinear contractions. U. P. B. Sci. Bull., Series A, 80(4), (2018), 21-28.}
\bibitem{SAM18}{\rm Shatanawi, W., Abodayeh, K., Mukheimer, A.:\;Some fixed point theorems in extended $b$-metric spaces. U. P. B. Sci. Bull., Series A, 80(4), 71-78 (2018).}
\bibitem{Ba89}{\rm Bakhtin, A.:\;The contraction mapping principle in almost metric spaces. Funct. Anal. 30, (1989), 26-37.}
\bibitem{KSA17}{\rm Kamran, T., Samreen, M., UL Ain, Q.:\;A generalization of $b$-metric space and some fixed point theorems. Mathematics 2017, 5, 19.}
\bibitem{MT94}{\rm Matthews, SG.:\;Partial metric topology: Ann. N.Y. Acad. Sci., 728(1), (1994), 183-197.}
\bibitem{Ba16}{\rm Bai, C.:\;Coupled fixed point theorems in $C^*$-algebra-valued $b$-metric spaces with application: Fixed Point Theory Appl. 2016, 2016:70.} 
\bibitem{BK15}{\rm Batul, S., Kamran, T.: $C^*$-valued contractive type mappings: Fixed Point Theory Appl. 2015, 2015:142.} 
\bibitem{EA18}{\rm Erden, M., Alaca, C.: $C^*$-algebra-valued $S$-metric spaces: Communications Series A1, 67(2), (2018),  165-177.}
\bibitem{KK18}{\rm Kalaivani, C., Kalpana, G: Fixed point theorems in $C^\ast$-algebra-valued $S$-metric spaces with some applications: U.P.B. Sci. Bull., Series A, 80(3), 2018.}
\bibitem{KS18}{\rm Kalpana, G., Sumaiya Tasneem, Z.:\;$C^*$-algebra-valued rectangular $b$-metric spaces and some fixed point theorems: Commun. Fac. Sci. Univ. Ank. Ser. A1 Math. Stat. 68(2), (2019), 2198-2208.}
\bibitem{KAS19}{\rm Kalpana, G., Sumaiya Tasneem, Z.:\;Common Fixed Point Theorems in $C^*$-algebra-valued Hexagonal $b$-Metric spaces: AIP Conf. Proc. 2095, 030012-1--030012-5, https://doi.org/10.1063/1.5097523.}
\bibitem{MG90}{\rm Murphy, G. J.: $C^*$-Algebras and Operator Theory:\;Academic Press, London 1990.}
\bibitem{QLZ16}{\rm Xin, QL., Jiang, LN., Ma, ZH.:\;Common fixed point theorems in $C^*$-algebra-valued metric spaces: J. Nonlinear Sci. Appl. 9, (2016), 4617-4627.}
\bibitem{SD16}{\rm Shehwar, D., Batul, S., Kamran, T., Ghiura, A.:\;Caristis fixed point theorem on $C^*$-algebra valued metric spaces: J. Nonlinear Sci. Appl., 9, (2016), 584-588.}
\bibitem{MZJ15}{\rm Ma, ZH., Jiang, LN.:\;$C^*$-algebra-valued $b$-metric spaces and related fixed point theorems: Fixed Pont Theory Appl. 2015, 2015:222.}
\bibitem{MAAS18}{\rm Mlaiki, N., Aydi H., Souayah, N., Abdeljawad:\;T. Controlled metric type spaces and the related contraction principle: Mathematics 2018, 6,\;194.}
\bibitem{ZLH14}{\rm Ma, ZH., Jiang, LN., Sun, H.:\;$C^*$-algebra-valued metric spaces and related fixed point theorems: Fixed Point Theory Appl. 2014, 2014:206.}
\end{thebibliography}
\end{document}